# Nonlinear Stabilization under Sampled and Delayed Measurements, and with Inputs Subject to Delay and Zero-Order Hold

Iasson Karafyllis[1] and Miroslav Krstic[2]

**Abstract**
Sampling arises simultaneously with input and output delays in networked control systems. When the delay is left uncompensated, the sampling period is generally required to be sufficiently small, the delay sufficiently short, and, for nonlinear systems, only semiglobal practical stability is generally achieved. For example, global stabilization of strict-feedforward systems under sampled measurements, sampled-data stabilization of the nonholonomic unicycle with arbitrarily sparse sampling, and sampled-data stabilization of LTI systems over networks with long delays, are open problems. In this paper we present two general results that address these example problems as special cases. First, we present global asymptotic stabilizers for forward complete systems under arbitrarily long input and output delays, with arbitrarily long sampling periods, and with continuous application of the control input. Second, we consider systems with sampled measurements and with control applied through a zero-order hold, under the assumption that the system is stabilizable under sampled-data feedback for some sampling period, and then construct sampled-data feedback laws that achieve global asymptotic stabilization under arbitrarily long input and measurement delays. All the results employ "nominal" feedback laws designed for the continuous-time systems in the absence of delays, combined with "predictor-based" compensation of delays and the effect of sampling.

**Keywords:** feedback stabilization, time-delay systems, sampled-data systems, nonlinear control.

## 1. Introduction

*Motivation*. Sampling arises simultaneously with input and output delays in many control problems, most notably in control over networks. In the absence of delays, in sampled-data control of nonlinear systems semiglobal practical stability is generally guaranteed [5,27,28,29], with the desired region of attraction achieved by sufficiently fast sampling. Alternatively, global results are achieved under restrictive conditions on the structure of the system [4,7,11,12,14,31].

[1] Dept. of Environmental Eng., Technical University of Crete, 73100, Chania, Greece, email: ikarafyl@enveng.tuc.gr.

[2] Dept. of Mechanical and Aerospace Eng., University of California, San Diego, La Jolla, CA 92093-0411, U.S.A., email: krstic@ucsd.edu.

On the other hand, in purely continuous-time nonlinear control, input delays of arbitrary length can be compensated [15,19,20] but no sampled-data extensions of such results are available. Simultaneous consideration to sampling and delays (either physical or sampling-induced) is given in the literature on control of linear and nonlinear systems over networks [2,3,6,26,30,31,34,35,37], but all available results rely on delay-dependent conditions for the existence of stabilizing feedback.

Despite the remarkable accomplishments in the fields of sampled-data, networked, and nonlinear delay systems, the following example problems remain open: global stabilization of strict-feedforward systems under sampled measurements and continuous control, sampled-data stabilization of the nonholonomic unicycle with inputs applied via zero-order hold and under arbitrarily sparse sampling, and sampled-data stabilization of LTI systems over networks with long delays.

In this paper we introduce two frameworks for solving such problems:

1. We present global asymptotic stabilizers for forward complete systems under arbitrarily long input and output delays, with arbitrarily long sampling periods, and with continuous application of the control input.
2. We consider systems with sampled measurements and with control applied through a zero-order hold, under the assumption that the system is stabilizable under sampled-data feedback for some sampling period, and then construct sampled-data feedback laws that achieve global asymptotic stabilization under arbitrarily long input and measurement delays.

In both frameworks we employ "nominal" feedback laws designed in the absence of delays, combined with "predictor-based" compensation of delays.

*Problem Statement.* As in [15,19,20,21,22,23,24,36,38], we consider systems with input delay,

$$\dot{x}(t) = f(x(t), u(t-\tau)) \tag{1.1}$$

where $x(t) = (x_1(t),...,x_n(t))' \in \Re^n$, $u(t) \in \Re^m$, $f : \Re^n \times \Re^m \to \Re^n$ is a locally Lipschitz mapping with $f(0,0) = 0$ and $\tau \geq 0$ is a constant. In [15,19,20,21,38], the feedback design problem for system (1.1) is addressed by assuming a feedback stabilizer $u = k(x)$ for system (1.1) with no delay, i.e. (1.1) with $\tau = 0$, or

$$\dot{x}(t) = f(x(t), u(t)) \tag{1.2}$$

and applying a delay compensator (predictor) methodology based on the knowledge of the delay. In this paper, we incorporate also a consideration of measurement delay, namely, we address the problem of stabilization of (1.1) with output

$$y(t) = x(t-r) \in \Re^n \tag{1.3}$$

where $r \geq 0$ is a constant, i.e., we consider delayed measurements. The motivation for a simultaneous consideration of input and measurement delays is that in many chemical process control problems the measurement delay of concentrations of chemical species can be large.

We also assume that the output is available at discrete time instants $\tau_i$ (the sampling times) with $\tau_{i+1} - \tau_i = T > 0$, where $T > 0$ is the sampling period. Very few papers have studied this problem (an exception is [8] where input and measurement delays are considered for linear systems but the measurement is not sampled).

The problem of stabilization of (1.1) with output given by (1.3) is intimately related to the stabilization of system (1.1) alone. To see this, notice that the output $y(t)$ of (1.1), (1.3) satisfies the following system of differential equations for all $t \geq r$:

$$\dot{y}(t) = f(y(t), u(t-r-\tau))$$

Consider the comparison between two problems described by the same differential equations: the problem of stabilization of (1.1) with input delay $r > 0$ and no measurement delay (i.e., $\dot{x}(t) = f(x(t), u(t-r))$ for all $t \geq 0$) and the problem of stabilization of (1.1), (1.3) with no input delay and measurement delay $r > 0$ (i.e., $\dot{y}(t) = f(y(t), u(t-r))$ for all $t \geq r$). The two problems are not identical: in the first stabilization problem the applied input values for $t \in [0, r]$ are given (as initial conditions), while in the second stabilization problem the applied input values for $t \in [0, r]$ must be computed based on an arbitrary initial condition $x(\theta) = x_0(\theta)$, $\theta \in [-r, 0]$ (irrespective of the current value of the state). Therefore, serious technical issues concerning the existence of the solution for $t \in [0, r]$ arise for the second stabilization problem (see Remark 2.2(b) below).

*Results of the paper.* We establish two general results:

1. A solution for the stabilization of (1.1) with output given by (1.3) under the assumption that system (1.2) is globally stabilizable and forward complete and *the input can be continuously adjusted* (Theorem 2.1). The proposed dynamic sampled-data controller uses values of the output (1.3) at the discrete time instants $\tau_i = t_0 + iT, i \in Z^+$, where $T > 0$ is the sampling period and $t_0 \geq 0$ is the initial time. This justifies the term "sampled-data". No restrictions for the values of the delays $r, \tau \geq 0$ or the sampling period $T > 0$ are imposed. In general, we show that there is no need for continuous measurements for global asymptotic stabilization of any stabilizable forward complete system with arbitrary input and output delays.

2. A solution for the stabilization of (1.1) with output given by (1.3) under the assumption that system (1.2) is globally stabilizable and forward complete and *the control action is implemented with zero order hold* (Theorem 3.2). Again, the proposed sampled-data controller uses values of the output (1.3) at the discrete time instants $\tau_i = t_0 + iT, i \in Z^+$, where $T > 0$ is the sampling period and $t_0 \geq 0$ is the initial time. In this case, we can solve the stabilization problem for systems with both delayed inputs and measurements provided that the user chooses the sampling period as the ratio of the input delay and any integer.

Our delay compensation methodology guarantees that any controller (continuous or sampled-data) designed for the delay-free case can be used for the regulation of the delayed system with input/measurement delays and sampled measurements. For example, all sampled-data feedback designs proposed in [4,5,11,14,27,28,29,31] which guarantee global stabilization can be exploited for the stabilization of a delayed system with input/measurement delays, sampled measurements and input applied with zero order hold.

The results are applied to
- the Linear Time Invariant (LTI) case, where $f(x,u) = Ax + Bu$, $A \in \Re^{n \times n}, B \in \Re^{n \times m}$. This case has been recently studied extensively in the context of linear Networked Control Systems, where various delays arise. Delay-dependent and/or sampling period-dependent sufficient conditions for the stabilization of Networked Control Systems have been proposed in the literature [2,3,6,26,30,31,34,35,37]. Here, we propose a linear delay compensator that guarantees exponential stability of the closed-loop system under the mild restriction that the user chooses the sampling period as the ratio of the input delay and any integer, with no additional restrictions for the delays (Corollary 3.4). The compensator is designed based on the knowledge of linear feedback stabilizer for the delay-free case.
- strict-feedforward systems [18,20,33], which are studied in Examples in 2.4 and 3.8.
- the stabilization of the nonholonomic integrator

$$\dot{x}_1 = u_1 \quad , \quad \dot{x}_2 = x_1 u_2 \quad , \quad \dot{x}_3 = u_2 \tag{1.4}$$

with both delayed inputs and measurements. The problem was recently studied in [17] in the presence of delays and in [4,29] in the presence of sampling. Here, our proposed dynamic sampled-data controller is applied with no restrictions for the value of the delays or the size of the sampling period. The stabilization problem is solved for the case where the inputs can be continuously adjusted (Corollary 4.1), as well as for the case where the inputs are applied with zero order hold (Proposition 4.2).

*Organization of the paper.* In Section 2 the main results concerning the case of the continuously adjusted input are stated and many comments and explanations are provided. In Section 3 the main results concerning the case of input applied with zero order hold are provided. Special results are provided for the case of linear autonomous systems and for the case of nonlinear systems which are diffeomorphically equivalent to a chain of integrators. Section 4 is devoted to the application of the obtained results to the stabilization of a three-wheeled vehicle with two independent rear motorized wheels (the nonholonomic integrator). Finally, in Section 5 we present the concluding remarks of the present work. The Appendix contains the proofs of certain results.

*Notation.* Throughout the paper we adopt the following notation:

* For a vector $x \in \Re^n$ we denote by $|x|$ its usual Euclidean norm, by $x'$ its transpose. For a real matrix $A \in \Re^{n \times m}$, $A' \in \Re^{m \times n}$ denotes its transpose and $|A| := \sup\{|Ax|; x \in \Re^n, |x|=1\}$ is its induced norm. $I \in \Re^{n \times n}$ denotes the identity matrix.

* $\Re^+$ denotes the set of non-negative real numbers. $Z^+$ denotes the set of non-negative integers. For every $t \geq 0$, $[t]$ denotes the integer part of $t \geq 0$, i.e., the largest integer being less or equal to $t \geq 0$.

* For the definition of the class of functions $KL$, see [16].

* By $C^j(A)$ ($C^j(A;\Omega)$), where $j \geq 0$ is a non-negative integer, we denote the class of functions (taking values in $\Omega$) that have continuous derivatives of order $j$ on $A$.

* Let $x:[a-r,b] \to \Re^n$ with $b > a \geq 0$ and $r \geq 0$. By $T_r(t)x$ we denote the "history" of $x$ from $t-r$ to $t$, i.e., $(T_r(t)x)(\theta) := x(t+\theta); \theta \in [-r,0]$, for $t \in [a,b]$. By $\breve{T}_r(t)x$ we denote the "open history" of $x$ from $t-r$ to $t$, i.e., $(\breve{T}_r(t)x)(\theta) := x(t+\theta); \theta \in [-r,0)$, for $t \in [a,b]$.

* Let $I \subseteq \Re^+ := [0,+\infty)$ be an interval. By $L^\infty(I;U)$ ($L^\infty_{loc}(I;U)$) we denote the space of measurable and (locally) bounded functions $u(\cdot)$ defined on $I$ and taking values in $U \subseteq \Re^m$. Notice that we do not identify functions in $L^\infty(I;U)$ which differ on a measure zero set. For $x \in L^\infty([-r,0];\Re^n)$ or $x \in L^\infty([-r,0);\Re^n)$ we define $\|x\|_r := \sup_{\theta \in [-r,0]} |x(\theta)|$ or $\|x\|_r := \sup_{\theta \in [-r,0)} |x(\theta)|$. Notice that $\sup_{\theta \in [-r,0]} |x(\theta)|$ is not

the essential supremum but the actual supremum and that is why the quantities $\sup_{\theta \in [-r,0]} |x(\theta)|$ and $\sup_{\theta \in [-r,0)} |x(\theta)|$ do not coincide in general. We will also use the notation $M_U$ for the space of measurable and locally bounded functions $u : \Re^+ \to U$.

* We say that a system of the form (1.2) is forward complete if for every $x_0 \in \Re^n$, $u \in M_U$ the solution $x(t)$ of (1.2) with initial condition $x(0) = x_0 \in \Re^n$ corresponding to input $u \in M_U$ exists for all $t \geq 0$.

Throughout the paper we adopt the convention $L^\infty([-r,0];\Re^n) = \Re^n$ and $C^0([-r,0];\Re^n) = \Re^n$ for $r = 0$. Finally, for reader's convenience, we mention the following fact, which is a direct consequence of Lemma 2.2 in [1] and Lemma 3.2 in [10]. The fact is used extensively throughout the paper.

**FACT:** *Suppose that system (1.2) is forward complete. Then for every $x_0 \in \Re^n$, $u \in L^\infty_{loc}([-\tau,+\infty);\Re^m)$ the solution $x(t)$ of (1.1) with initial condition $x(0) = x_0 \in \Re^n$ corresponding to input $u \in L^\infty_{loc}([-\tau,+\infty);\Re^m)$ exists for all $t \geq 0$. Moreover, for every $T > 0$ there exists a function $a \in K_\infty$ such that for every $x_0 \in \Re^n$, $u \in L^\infty_{loc}([-\tau,+\infty);\Re^m)$ the solution $x(t)$ of (1.1) with initial condition $x(0) = x_0 \in \Re^n$ corresponding to input $u \in L^\infty_{loc}([-\tau,+\infty);\Re^m)$ satisfies*

$$|x(t)| \leq a\left(|x_0| + \sup_{-\tau \leq s \leq t-\tau} |u(s)|\right), \text{ for all } t \in (0,T].$$

## 2. Dynamic Sampled-Data Feedback for Continuously Adjusted Input

We start by presenting the assumptions for system (1.2). Our first assumption concerning system (1.2) is forward completeness.

**Hypothesis (H1):** *System (1.2) is forward complete.*

Assumption (H1) guarantees that system (1.1) is forward complete as well: for every $x_0 \in \Re^n$, $u \in L^\infty_{loc}([-\tau,+\infty);\Re^m)$ the solution $x(t)$ of (1.1) with initial condition $x(0) = x_0 \in \Re^n$ corresponding to input $u \in L^\infty_{loc}([-\tau,+\infty);\Re^m)$ exists for all $t \geq 0$. Therefore, we are in a position to define the "predictor" mapping $\Phi : \Re^n \times L^\infty([-r-\tau,0);\Re^m) \to \Re^n$ for all $r,\tau \geq 0$ with $r + \tau > 0$ in the following way:

"for every $x_0 \in \Re^n$, $u \in L^\infty\big([-r-\tau,0);\Re^m\big)$ the solution $x(t)$ of (1.1) with initial condition $x(-r) = x_0$ corresponding to input $u \in L^\infty\big([-r-\tau,0);\Re^m\big)$ satisfies $x(\tau) = \Phi(x_0,u)$"

By virtue of the Fact, we can guarantee the existence of $a \in K_\infty$ such that

$$|\Phi(x,u)| \leq a\big(|x| + \|u\|_{r+\tau}\big), \text{ for all } (x,u) \in \Re^n \times L^\infty\big([-r-\tau,0);\Re^m\big) \tag{2.1}$$

Using (2.1) and the fact that $f : \Re^n \times \Re^m \to \Re^n$ is a locally Lipschitz mapping, we can guarantee the existence of a non-decreasing function $L : \Re^+ \to \Re^+$ such that

$$|\Phi(x,u) - \Phi(y,v)| \leq L\big(|x| + |y| + \|u\|_{r+\tau} + \|v\|_{r+\tau}\big)\big(|x-y| + \|u-v\|_{r+\tau}\big),$$

$$\text{for all } (x,u) \in \Re^n \times L^\infty\big([-r-\tau,0);\Re^m\big), (y,v) \in \Re^n \times L^\infty\big([-r-\tau,0);\Re^m\big) \tag{2.2}$$

We assume next that (1.2) is globally stabilizable.

**Hypothesis (H2) (continuously adjusted input):** *There exists $k \in C^1(\Re^+ \times \Re^n; \Re^m)$, $g \in K_\infty$ with*

$$|k(t,x)| \leq g(|x|), \text{ for all } (t,x) \in \Re^+ \times \Re^n \tag{2.3}$$

*such that $0 \in \Re^n$ is Uniformly Globally Asymptotically Stable for system (1.2) with $u = k(t,x)$, i.e., there exists a function $\sigma \in KL$ such that for every $(t_0,x_0) \in \Re^+ \times \Re^n$ the solution $x(t)$ of (1.2) with $u = k(t,x)$ and initial condition $x(t_0) = x_0 \in \Re^n$ satisfies the following inequality:*

$$|x(t)| \leq \sigma\big(|x_0|, t-t_0\big), \forall t \geq t_0 \tag{2.4}$$

Consider system (1.1) under hypotheses (H1), (H2) for system (1.2). Our proposed dynamic sampled-data feedback has states $(z(t), T_{r+\tau}(t)u) \in \Re^n \times L^\infty\big([-r-\tau,0];\Re^m\big)$ and inputs $y(t) \in \Re^n$ and for each $t_0 \geq 0$, $(z_0, u_0) \in \Re^n \times L^\infty\big([-r-\tau,0];\Re^m\big)$ the states are computed by the interconnection of two subsystems:

1) A sampled-data subsystem (see [10]) with inputs $(y(t), T_{r+\tau}(t)u) \in \Re^n \times L^\infty\big([-r-\tau,0];\Re^m\big)$:

$$\begin{aligned} \dot{z}(t) &= f(z(t), u(t)), \ t \in [\tau_i, \tau_{i+1}), i \in Z^+ \\ z(\tau_{i+1}) &= \Phi\big(y(\tau_{i+1}), \breve{T}_{r+\tau}(\tau_{i+1})u\big) \\ z(t_0) &= z_0 \in \Re^n \end{aligned} \tag{2.5}$$

where

$$\tau_i = t_0 + iT, \ i \in Z^+$$

are the sampling times and $T > 0$ is the sampling period. We stress that the proposed sampled-data dynamic controller uses only values of the output $y(t) = x(t-r) \in \Re^n$ at the discrete time instants $\tau_i = t_0 + iT$, where $i \in Z^+$.

2) A subsystem described by Functional Difference Equations (see [13]) with inputs $z(t) \in \Re^n$:

$$u(t) = k(t+\tau, z(t)), t > t_0$$
$$T_{r+\tau}(t_0)u = u_0 \in L^\infty([-r-\tau, 0]; \Re^m)$$
(2.6)

Our first main result is now stated.

**Theorem 2.1:** *Let $T > 0$, $r, \tau \geq 0$ with $r + \tau > 0$ and suppose that hypotheses (H1), (H2) hold for system (1.2). Then the closed-loop system (1.1), (1.3) (2.5), (2.6) is Uniformly Globally Asymptotically Stable, in the sense that there exists a function $\tilde{\sigma} \in KL$ such that for every $t_0 \geq 0$, $(x_0, z_0, u_0) \in C^0([-r,0]; \Re^n) \times \Re^n \times L^\infty([-r-\tau, 0]; \Re^m)$, the solution $(x(t), z(t), u(t)) \in \Re^n \times \Re^n \times \Re^m$ of the closed-loop system (2.5), (2.6), (1.3), (1.1) with initial condition $z(t_0) = z_0 \in \Re^n$, $T_{r+\tau}(t_0)u = u_0 \in L^\infty([-r-\tau, 0]; \Re^m)$, $T_r(t_0)x = x_0 \in C^0([-r,0]; \Re^n)$ satisfies the following inequality for all $t \geq t_0$:*

$$|z(t)| + \|T_r(t)x\|_r + \|T_{r+\tau}(t)u\|_{r+\tau} \leq \tilde{\sigma}(|z_0| + \|x_0\|_r + \|u_0\|_{r+\tau}, t - t_0)$$
(2.7)

Some remarks for the dynamic sampled-data feedback given by (2.5), (2.6) are in order before we proceed to the proof of Theorem 2.1.

**Remark 2.2:**
**(a)** The dynamic sampled-data controller (2.5), (2.6) is time-varying if $k$ is time-varying. If $k$ is $T$–periodic then the dynamic sampled-data controller (2.5), (2.6) is $T$–periodic too.
**(b)** Since the output $y(t)$ given by (1.3) of system (1.1) satisfies the system of differential equations $\dot{y}(t) = f(y(t), u(t-r-\tau))$, for all $t \geq t_0 + r$, where $t_0 \geq 0$ is the initial time, we can in principle apply the predictor-based delay compensation approach described in [20] (extended for time-varying feedback laws), which gives the static feedback law $u(t) = k(t+\tau, \Phi(y(t), \breve{T}_{r+\tau}(t)u))$. This is the inspiration for the construction of the sampled-data dynamic feedback (2.5), (2.6): for all $t = t_0 + iT$, where $i \in Z^+$, the value of $u(t)$ computed by (2.5), (2.6) is exactly $u(t) = k(t+\tau, \Phi(y(t), \breve{T}_{r+\tau}(t)u))$. However, there are technical problems with the application of the static feedback law $u(t) = k(t+\tau, \Phi(y(t), \breve{T}_{r+\tau}(t)u))$: given initial conditions $T_r(t_0)x = x_0 \in C^0([-r,0]; \Re^n)$, $T_{r+\tau}(t_0)u = u_0 \in L^\infty([-r-\tau, 0]; \Re^m)$, we cannot guarantee existence of a measurable and essentially bounded $u: [t_0 - r - \tau, t_0 + r] \to \Re^m$ satisfying the integral equation $u(t) = k(t+\tau, \Phi(x(t-r), \breve{T}_{r+\tau}(t)u))$ for all $t \in (t_0, t_0 + r]$. For the case $\tau = 0$, one sufficient condition for the existence of a solution of the integral equation is that the initial condition $T_r(t_0)x = x_0 \in C^0([-r,0];\Re^n)$, $T_r(t_0)u = u_0 \in L^\infty([-r,0];\Re^m)$ satisfies the equation $\dot{x}(t-r) = f(x(t-r),u(t-r))$ for all $t \in [t_0, t_0+r]$: in this case the solution is $u(t) = k(t,z(t))$ for $t \in (t_0, t_0+r]$, where $z(t)$ is the solution of the initial value problem $\dot{z}(t) = f(z(t),k(t,z(t)))$ with $z(t_0) = x(t_0)$. Other restrictive sufficient conditions for the existence of a solution of the integral equation can be obtained by using fixed point theory. The proof of Theorem 2.1 shows that this issue can be completely avoided for the dynamic sampled-data feedback (2.5), (2.6).

**(c)** For every initial condition the value of $u(t)$ computed by (2.5), (2.6) is exactly $u(t) = k(t+\tau, \Phi(y(t), \breve{T}_{r+\tau}(t)u))$ for all $t \geq t_0 + iT$ with $i \in Z^+$ satisfying $iT \geq r$, so our dynamic sampled-data feedback is based on the predictor principle.

**(d)** For the implementation of the controller (2.5), (2.6), we must know the "predictor" mapping $\Phi : \Re^n \times L^\infty([-r-\tau,0);\Re^m) \to \Re^n$. This mapping can be explicitly computed for

(i) Linear systems $\dot{x} = Ax + Bu$, with $x \in \Re^n, u \in \Re^m$. In this case (Corollary 3.4 below) the predictor mapping $\Phi : \Re^n \times L^\infty([-r-\tau,0);\Re^m) \to \Re^n$ is given by the explicit equation

$$\Phi(x,u) := \exp(A(\tau+r))x + \int_{-r-\tau}^{0} \exp(-Aw)Bu(w)dw.$$

(ii) Bilinear systems $\dot{x} = Ax + Bu + uCx$, with $x \in \Re^n, u \in \Re$ and $AC = CA$. In this case the predictor mapping $\Phi : \Re^n \times L^\infty([-r-\tau,0);\Re^m) \to \Re^n$ is given by the explicit equation

$$\Phi(x,u) := \exp(A(\tau+r))\exp\left(C\int_{-r-\tau}^{0} u(s)ds\right)x + \int_{-r-\tau}^{0} \exp(-Aw)\exp\left(C\int_{w}^{0} u(s)ds\right)Bu(w)dw.$$

(iii) Nonlinear systems of the following form:

$$\begin{aligned}\dot{x}_1 &= a_1(u)x_1 + f_1(u)\\ \dot{x}_2 &= a_2(u,x_1)x_2 + f_2(u,x_1)\\ &\vdots\\ \dot{x}_n &= a_n(u,x_1,...,x_{n-1})x_n + f_n(u,x_1,...,x_{n-1})\\ x &= (x_1,...,x_n)' \in \Re^n, u \in \Re^m\end{aligned}$$

where all mappings $a_i, f_i$ ($i = 1,...,n$) are locally Lipschitz. In this case the predictor mapping $\Phi : \Re^n \times L^\infty([-r-\tau,0);\Re^m) \to \Re^n$ can be constructed inductively. For example, for $n=1$ the predictor mapping is given by $\Phi(x,u) = \exp\left(\int_{-r-\tau}^{0} a_1(u(s))ds\right)x + \int_{-r-\tau}^{0} \exp\left(\int_{w}^{0} a_1(u(s))ds\right)f_1(u(w))dw$.

Example 2.4 below applies Theorem 2.1 to a three-dimensional nonlinear system of the above class. Moreover, the nonholonomic integrator (1.4) belongs to the above class and Theorem 2.1 can be applied (see Corollary 4.1).

(iv) Nonlinear systems $\dot{x} = f(x,u)$, for which there exists a global diffeomorphism $\Theta: \Re^n \to \Re^n$ such that the change of coordinates $z = \Theta(x)$ transforms the system to one of the above cases (Corollary 3.7 below).

For globally Lipschitz systems, one can utilize approximate "predictor" mappings $\Phi: \Re^n \times L^\infty([-r-\tau,0); \Re^m) \to \Re^n$ as shown in [15] under additional and more restrictive hypotheses.

**Proof of Theorem 2.1:** We start with the following claim, which we prove in the Appendix.

**Claim 1:** There exists a function $G \in K_\infty$ such that for every $t_0 \geq 0$ and $(x_0, z_0, u_0) \in C^0([-r,0]; \Re^n) \times \Re^n \times L^\infty([-r-\tau,0]; \Re^m)$ the solution $(x(t), z(t), u(t)) \in \Re^n \times \Re^n \times \Re^m$ of the closed-loop system (2.5), (2.6), (1.3), (1.1) with initial condition $z(t_0) = z_0 \in \Re^n$, $T_{r+\tau}(t_0)u = u_0 \in L^\infty([-r-\tau,0]; \Re^m)$, $T_r(t_0)x = x_0 \in C^0([-r,0]; \Re^n)$ exists for all $t \in [t_0, t_0 + T]$ and satisfies

$$\|T_{r+\tau}(t)u\|_{r+\tau} + |z(t)| + \|T_r(t)x\|_r \leq G(|z_0| + \|x_0\|_r + \|u_0\|_{r+\tau}), \text{ for all } t \in [t_0, t_0 + T] \quad (2.8)$$

By virtue of induction and Claim 1, the following claim holds.

**Claim 2:** There exists a function $G \in K_\infty$ such that for every $t_0 \geq 0$, $p \in Z^+$, $p \geq 1$ and $(x_0, z_0, u_0) \in C^0([-r,0]; \Re^n) \times \Re^n \times L^\infty([-r-\tau,0]; \Re^m)$ the solution $(x(t), z(t), u(t)) \in \Re^n \times \Re^n \times \Re^m$ of the closed-loop system (2.5), (2.6), (1.3), (1.1) with initial condition $z(t_0) = z_0 \in \Re^n$, $T_{r+\tau}(t_0)u = u_0 \in L^\infty([-r-\tau,0]; \Re^m)$, $T_r(t_0)x = x_0 \in C^0([-r,0]; \Re^n)$ exists for all $t \in [t_0, t_0 + pT]$ and satisfies

$$\|T_{r+\tau}(t)u\|_{r+\tau} + |z(t)| + \|T_r(t)x\|_r \leq G^{(p)}(|z_0| + \|x_0\|_r + \|u_0\|_{r+\tau}), \text{ for all } t \in [t_0, t_0 + pT] \quad (2.9)$$

where $G^{(p)}(s) := \underbrace{G \circ \ldots \circ G}_{p \text{ times}}$.

We notice that for all $t \geq t_0 + iT$ with $i \in Z^+$ satisfying $iT \geq r$ the solution $(x(t), z(t), u(t)) \in \Re^n \times \Re^n \times \Re^m$ of the closed-loop system (2.5), (2.6), (1.3), (1.1) satisfies:

$$z(t) = x(t + \tau) \quad (2.10)$$

Consequently, for all $t \geq t_0 + iT + \tau$ with $i \in Z^+$ satisfying $iT \geq r$ it holds that:

$$u(t - \tau) = k(t, x(t)) \quad (2.11)$$

Hypothesis (H2) in conjunction with inequality (2.4) and equation (2.11) implies that the following inequality holds:

$$|x(t)| \leq \sigma(|x(t_0 + iT + \tau)|, t - t_0 - iT - \tau), \forall t \geq t_0 + iT + \tau \quad (2.12)$$

Define $p = \left[\frac{r}{T}\right] + \left[\frac{\tau}{T}\right] + 2$. Using (2.9), (2.10) and (2.12), it follows that the following inequality holds:

$$|x(t)| + |z(t)| \leq \sigma\left(G^{(p)}\left(\|T_{r+\tau}(t_0)u\|_{r+\tau} + |z(t_0)| + \|T_r(t_0)x\|_r\right), t - t_0 - pT\right), \quad \forall t \geq t_0 + pT \tag{2.13}$$

Define $\tilde{\sigma}(s,t) := \sigma\left(G^{(p)}(s), t - pT - r\right) + g\left(\sigma\left(G^{(p)}(s), t - pT - r\right)\right)$ for all $t \geq pT + r$ and $\tilde{\sigma}(s,t) := G^{(p)}(s) + g\left(G^{(p)}(s)\right)$ for all $t \in [0, pT + r)$. Using (2.3), (2.10), (2.11) and (2.13) we can conclude that (2.7) holds. The proof is complete. ◁

**Remark 2.3:** The proof of Theorem 2.1 shows that if $k \in C^1(\Re^+ \times \Re^n; \Re^m)$ with $k(t,0) = 0$ for all $t \geq 0$, is a non-uniform in time stabilizer for system (see [9]) then we can prove that the closed-loop system (2.5), (2.6), (1.3), (1.1) is non-uniformly in time Globally Asymptotically Stable, i.e., there exist functions $\tilde{\sigma} \in KL$ and a positive continuous function $\beta : \Re^+ \to (0, +\infty)$ such that for every $t_0 \geq 0$, $(x_0, z_0, u_0) \in C^0([-r,0];\Re^n) \times \Re^n \times L^\infty([-r-\tau,0];\Re^m)$, the solution $(x(t), z(t), u(t)) \in \Re^n \times \Re^n \times \Re^m$ of the closed-loop system (2.5), (2.6), (1.3), (1.1) with initial condition $z(t_0) = z_0 \in \Re^n$, $T_{r+\tau}(t_0)u = u_0 \in L^\infty([-r-\tau,0];\Re^m)$, $T_r(t_0)x = x_0 \in C^0([-r,0];\Re^n)$ satisfies the following inequality for all $t \geq t_0$:

$$|z(t)| + \|T_r(t)x\|_r + \|T_{r+\tau}(t)u\|_{r+\tau} \leq \tilde{\sigma}\left(\beta(t_0)\left(|z_0| + \|x_0\|_r + \|u_0\|_{r+\tau}\right), t - t_0\right) \tag{2.14}$$

In this case, there is no need to assume that inequality (2.3) holds.

We next present an example which shows how the obtained results can be applied to feedforward nonlinear systems.

**Example 2.4 (Control of strict-feedforward systems with arbitrarily sparse sampling):**
Consider the following example taken from [20]:

$$\dot{x}_1(t) = x_2(t) + x_3^2(t), \quad \dot{x}_2(t) = x_3(t) + x_3(t)u(t-\tau), \quad \dot{x}_3(t) = u(t-\tau) \tag{2.15}$$
$$x(t) = (x_1(t), x_2(t), x_3(t))' \in \Re^3, u(t) \in \Re$$

Here, we consider the stabilization problem for (2.15) with output given by (1.3) available only at the discrete time instants $\tau_i$ (the sampling times) with $\tau_{i+1} - \tau_i = T > 0$, where $T > 0$ is the sampling period. Hypothesis (H1) holds for system (2.15) and the predictor mapping can be explicitly expressed by the equations:

$$\Phi(x,u) := \left[\phi_1(x,u), \quad x_2 + (\tau+r)x_3 + x_3 \int_{-r-\tau}^{0} u(s)ds + \int_{-r-\tau}^{0} (1+u(s)) \int_{-r-\tau}^{s} u(q)dq\, ds, \quad x_3 + \int_{-r-\tau}^{0} u(s)ds\right]' \tag{2.16}$$

where

$$\phi_1(x,u) = x_1 + (\tau+r)x_2 + (\tau+r)x_3^2 + \frac{1}{2}(\tau+r)^2 x_3 + 3x_3 \int_{-r-\tau}^{0}\int_{-r-\tau}^{s} u(q)dq\, ds$$
$$+ \int_{-r-\tau}^{0}\int_{-r-\tau}^{s}(1+u(w))\int_{-r-\tau}^{w}u(q)dq\, dw\, ds + \int_{-r-\tau}^{0}\left(\int_{-r-\tau}^{s}u(q)dq\right)^2 ds \quad (2.17)$$

Moreover, hypothesis (H2) holds as well with the smooth, time-independent feedback law:

$$k(x) := -x_1 - 3x_2 - \frac{3}{8}x_2^2 + \frac{3}{4}x_3\left(-4 - x_1 - 2x_2 + \frac{1}{2}x_3 + \frac{1}{2}x_2 x_3 + \frac{5}{8}x_3^2 - \frac{1}{4}x_3^3 - \frac{3}{8}\left(x_2 - \frac{1}{2}x_3^2\right)^2\right) \quad (2.18)$$

It follows from Theorem 2.1 that the dynamic sampled-data controller $u(t) = k(z(t))$ with

$$\dot{z}_1(t) = z_2(t) + z_3^2(t),\ \dot{z}_2(t) = z_3(t) + z_3(t)u(t),\ \dot{z}_3(t) = u(t)$$
$$z(t) = (z_1(t), z_2(t), z_3(t))' \in \mathfrak{R}^3, \text{ for } t \in [\tau_i, \tau_{i+1}) \quad (2.19)$$

and

$$z(\tau_{i+1}) = \Phi(y(\tau_{i+1}), \breve{T}_{r+\tau}(\tau_{i+1})u)\ , i \in Z^+ \quad (2.20)$$

where $\Phi : \mathfrak{R}^3 \times L^\infty([-r-\tau,0);\mathfrak{R}^m) \to \mathfrak{R}^3$ is defined by (2.16), (2.17) and $k : \mathfrak{R}^3 \to \mathfrak{R}$ is defined by (2.18), guarantees global asymptotic stability for system (2.15). The reader should notice that the dynamic sampled-data controller (2.19), (2.20), (2.21) can still be used even if no delays are present but the state is available only at the discrete time instants $\tau_i$ (the sampling times) with $\tau_{i+1} - \tau_i = T > 0$, where $T > 0$ is the sampling period. Hence, in this section we have provided, as a special case, the first solution to the problem of global asymptotic stabilization of strict feedforward systems with arbitrarily sparse in time sampling of the state and with continuous control. ◁

## 3. Sampled-Data Feedback for Input Applied with Zero Order Hold

### 3.1 General Design

This section is devoted to the case where the input is applied with zero order hold. In this section we assume that (1.2) is globally stabilizable with feedback applied with zero order hold. This is very different from hypothesis (H2) in the previous section.

**Hypothesis (H3) (input applied with zero order hold):** *There exists* $k : \mathfrak{R}^n \to \mathfrak{R}^m$, $g \in K_\infty$, $T > 0$ *such that*

$$|k(x)| \le g(|x|), \text{ for all } x \in \mathfrak{R}^n \quad (3.1)$$

and such that $0 \in \Re^n$ is Uniformly Globally Asymptotically Stable for the sampled-data system

$$\begin{aligned}
&\dot{x}(t) = f(x(t), k(x(\tau_i))) \, , \, t \in [\tau_i, \tau_{i+1}) \\
&x(\tau_{i+1}) = \lim_{t \to \tau_{i+1}^-} x(t) \\
&\tau_{i+1} = \tau_i + T \\
&\tau_0 = 0 \geq 0 \, , \, x(0) = x_0 \in \Re^n
\end{aligned} \quad (3.2)$$

in the sense that there exists a function $\sigma \in KL$ such that for every $x_0 \in \Re^n$ the solution $x(t)$ of (3.2) with initial condition $x(0) = x_0 \in \Re^n$ satisfies inequality (2.4) with $t_0 = 0$ for all $t \geq 0$.

**Remark 3.1:** Hypothesis (H3) seems like a restrictive hypothesis, because it demands global stabilizability by means of sampled-data feedback with positive sampling rate. However, hypothesis (H3) can be satisfied for:

(i) Linear stabilizable systems, where $f(x,u) = Ax + Bu$, $A \in \Re^{n \times n}, B \in \Re^{n \times m}$ (see Corollary 3.4 and Remark 3.5 below),

(ii) Nonlinear systems of the form $\dot{x} = f(x) + g(x)u$, $x \in \Re^n, u \in \Re$, where the vector field $f : \Re^n \to \Re^n$ is globally Lipschitz and the vector field $g : \Re^n \to \Re^n$ is locally Lipschitz and bounded, which can be stabilized by a globally Lipschitz feedback law $u = k(x)$ (see [7]).

(iii) Nonlinear systems of the form $\dot{x}_i = f_i(x,u) + g_i(x,u)x_{i+1}$ for $i = 1,...,n-1$ and $\dot{x}_n = f_n(x,u) + g_n(x,u)u$, where the drift terms $f_i(x,u)$ ($i = 1,...,n$) satisfy the linear growth conditions $|f_i(x)| \leq L|x_1| + ... + L|x_i|$ ($i = 1,...,n$) for certain constant $L \geq 0$ and there exist constants $b \geq a > 0$ such that $a \leq g_i(x,u) \leq b$ for all $i = 1,...,n$, $x \in \Re^n, u \in \Re$ (see [12]).

(iv) Asymptotically controllable homogeneous systems with positive minimal power and zero degree (see [4]).

(v) Systems satisfying the reachability hypotheses of Theorem 3.1 in [14], or hypotheses (A1), (A2), (A3) in Section 4 of [11],

(vi) Nonlinear systems $\dot{x} = f(x,u)$, for which there exists a global diffeomorphism $\Theta : \Re^n \to \Re^n$ such that the change of coordinates $z = \Theta(x)$ transforms the system to one of the above cases.

Consider system (1.1) under hypotheses (H1), (H3) for system (1.2). In this case we propose a feedback law that is simply a composition of the feedback stabilizer and the delay compensator:

$$u(t) = k\left(\Phi\left(y(\tau_i), \breve{T}_{r+\tau}(\tau_i)u\right)\right) \, , \, t \in [\tau_i, \tau_{i+1}) \quad (3.3)$$

where $\tau_i = iT, i \in Z^+$ are the sampling times and $\Phi : \Re^n \times L^\infty\left([-r-\tau,0); \Re^m\right) \to \Re^n$ is the predictor mapping involved in (2.1), (2.2). The control action is applied with zero order hold, i.e., it is constant on $[\tau_i, \tau_{i+1})$; however the control action affecting system (1.1) remains constant on the interval $[\tau_i + \tau, \tau_{i+1} + \tau)$.

Our main result is stated next.

**Theorem 3.2:** *Let $T > 0$, $r, \tau \geq 0$ with $r + \tau > 0$ and suppose that there exists $l \in Z^+$ such that $\tau = lT$. Moreover, suppose that hypotheses (H1), (H2) hold for system (1.2). Then the closed-loop system (1.1) with (3.3), i.e., the following sampled-data system*

$$\begin{aligned} &\dot{x}(t) = f(x(t), u(t-\tau)) \\ &u(t) = k\big(\Phi(x(\tau_i - r), \breve{T}_{r+\tau}(\tau_i)u)\big), \ t \in [\tau_i, \tau_{i+1}), i \in Z^+ \\ &\tau_{i+1} = \tau_i + T, \tau_0 = 0 \end{aligned} \quad (3.4)$$

*is Uniformly Globally Asymptotically Stable, in the sense that there exists a function $\tilde{\sigma} \in KL$ such that for every $(x_0, u_0) \in C^0([-r, 0]; \Re^n) \times L^\infty([-r-\tau, 0); \Re^m)$, the solution $(x(t), u(t)) \in \Re^n \times \Re^m$ of system (3.4) with initial condition $\breve{T}_{r+\tau}(0)u = u_0 \in L^\infty([-r-\tau, 0); \Re^m)$, $T_r(0)x = x_0 \in C^0([-r, 0]; \Re^n)$ satisfies the following inequality for all $t \geq 0$:*

$$\|T_r(t)x\|_r + \|\breve{T}_{r+\tau}(t)u\|_{r+\tau} \leq \tilde{\sigma}\big(\|x_0\|_r + \|u_0\|_{r+\tau}, t\big) \quad (3.5)$$

*Finally, if system (3.2) satisfies the dead-beat property of order $jT$, where $j \in Z^+$ is positive, i.e., for all $x_0 \in \Re^n$ the solution $x(t)$ of (3.2) with initial condition $x(0) = x_0 \in \Re^n$ satisfies $x(t) = 0$ for all $t \geq jT$ then system (3.4) satisfies the dead-beat property of order $\left(j + l + \left[\frac{r}{T}\right] + 1\right)T$, where $\left[\frac{r}{T}\right]$ is the integer part of $\frac{r}{T}$, i.e., for every $(x_0, u_0) \in C^0([-r, 0]; \Re^n) \times L^\infty([-r-\tau, 0); \Re^m)$, the solution $(x(t), u(t)) \in \Re^n \times \Re^m$ of system (3.4) with initial condition $\breve{T}_{r+\tau}(0)u = u_0 \in L^\infty([-r-\tau, 0); \Re^m)$, $T_r(0)x = x_0 \in C^0([-r, 0]; \Re^n)$ satisfies $x(t) = 0$ for all $t \geq \left(j + l + \left[\frac{r}{T}\right] + 1\right)T$.*

**Remark 3.3:**

**(a)** If we denote $T_1 \geq 0$ to be the delay in receiving the measured data, $T_2 \geq 0$ the computation time for the quantity $v = k(\Phi(x, u))$, where $(x, u) \in \Re^n \times L^\infty([-r-\tau, 0); \Re^m)$, and $T_3 \geq 0$ the time for the data to reach the actuator, then one should notice that $r = T_1 + T_2$ and $\tau = T_3$.

**(b)** In practice, when $r, \tau \geq 0$ are given, the control practitioner should look for sampled-data stabilizers satisfying hypothesis (H3) for certain sampling period $T > 0$ with $\tau = lT$ for some $l \in Z^+$. Therefore, the value of the sampling period is determined after the estimation of the input delay. If $\tau = 0$ then one can choose any sampling period $T > 0$ (the condition $\tau = lT$ holds with $l = 0$). However, if $\tau > 0$ then the sampling period is constrained to be less or equal to $\tau > 0$. This may seem to be impractical for the cases where $\tau > 0$ is small. However, in practice, a delay that is smaller than a reasonable sampling period would be typically ignored, or the control engineer can induce an input delay of magnitude equal to the sampling period (for example, by delaying the transmission of the control action).

**Proof of Theorem 3.2:** Using the Fact, we can guarantee the existence of $b \in K_\infty$ such that for every $(x_0, u_0) \in C^0([-r, 0]; \Re^n) \times L^\infty([-r-\tau, 0); \Re^m)$ the solution $x(t) \in \Re^n$ of (1.1) with initial condition $\breve{T}_{r+\tau}(0)u = u_0 \in L^\infty([-r-\tau, 0); \Re^m)$, $T_r(0)x = x_0 \in C^0([-r, 0]; \Re^n)$ exists for all $t \in [0, \tau]$ and satisfies

$$|x(t)| \leq b\left(|x(0)| + \|u_0\|_{r+\tau}\right), \text{ for all } t \in [0, \tau] \tag{3.6}$$

If $\tau > 0$ then the input $u(t)$ takes exactly $l$ values on the interval $t \in [0, \tau)$ with input values $u_i = u(t)$ for $t \in [(i-1)T, iT)$, $i = 1, \ldots, l$. We also set $u_{l+1} = u(lT) = u(\tau)$. Using (3.1) and (2.1), we obtain the following estimates:

$$|u_1| \leq g\left(a\left(|x(-r)| + \|u_0\|_{r+\tau}\right)\right)$$

$$|u_i| \leq g\left(a\left(|x((i-1)T - r)| + \|u_0\|_{r+\tau} + \max_{1 \leq q \leq i-1}|u_q|\right)\right), \text{ for } i = 2, \ldots, l+1$$

Using the above inequalities, the trivial inequality $|x(-r)| \leq \|x_0\|_r$ and the inequalities $|x((i-1)T - r)| \leq \|x_0\|_r + b(\|x_0\|_r + \|u_0\|_{r+\tau})$ for $i = 2, \ldots, l+1$ (which are direct consequences of (3.6)), we can construct a function $h \in K_\infty$ such that

$$|u_i| \leq h\left(\|x_0\|_r + \|u_0\|_{r+\tau}\right), \text{ for } i = 1, \ldots, l+1 \tag{3.7}$$

Thus we may conclude that there exists $H \in K_\infty$ such that

$$\|T_r(t)x\|_r + \|\breve{T}_{r+\tau}(t)u\|_{r+\tau} \leq H\left(\|x_0\|_r + \|u_0\|_{r+\tau}\right), \text{ for all } t \in [0, \tau] \tag{3.8}$$

Inequality (3.8) holds trivially for the case $\tau = 0$.

We next continue with the following claim. Its proof is provided in the Appendix.

**Claim 3:** There exists a function $G \in K_\infty$ such that the solution $(x(t), u(t)) \in \Re^n \times \Re^m$ of system (3.4) exists for all $t \in [\tau, T+\tau]$ and satisfies

$$\|\breve{T}_{r+\tau}(t)u\|_{r+\tau} + \|T_r(t)x\|_r \leq G\left(\|\breve{T}_{r+\tau}(\tau)u\|_{r+\tau} + \|T_r(\tau)x\|_r\right), \text{ for all } t \in [\tau, T+\tau] \tag{3.9}$$

By virtue of induction and Claim 3, the following claim holds.

**Claim 4:** There exists a function $G \in K_\infty$ such that for every $p \in Z^+$, $p \geq 1$ the solution $(x(t), u(t)) \in \Re^n \times \Re^m$ of system (3.4) exists for all $t \in [\tau, \tau + pT]$ and satisfies

$$\left\| \breve{T}_{r+\tau}(t)u \right\|_{r+\tau} + \left\| T_r(t)x \right\|_r \leq G^{(p)}\left( \left\| \breve{T}_{r+\tau}(\tau)u \right\|_{r+\tau} + \left\| T_r(\tau)x \right\|_r \right), \text{ for all } t \in [\tau, pT + \tau] \quad (3.10)$$

where $G^{(p)}(s) := \underbrace{G \circ \ldots \circ G}_{p \text{ times}}$.

We notice that for all $i \in Z^+$ with $iT \geq r$ the solution $(x(t), u(t)) \in \Re^n \times \Re^m$ of system (3.4) satisfies:

$$u(t - \tau) = k(x(iT + \tau)), \quad \forall t \in [iT + \tau, (i+1)T + \tau) \quad (3.11)$$

Hypothesis (H3) in conjunction with inequality (2.4) with $t_0 = 0$ and equation (3.11) implies that the following inequality holds for all $i \in Z^+$ with $iT \geq r$:

$$|x(t)| \leq \sigma(|x(iT + \tau)|, t - iT - \tau), \quad \forall t \geq iT + \tau \quad (3.12)$$

Define $p = [r/T] + 1$. Using (3.10) and (3.12), the following inequality holds:

$$|x(t)| \leq \sigma\left( G^{(p)}\left( \left\| T_{r+\tau}(\tau)u \right\|_{r+\tau} + \left\| \breve{T}_r(\tau)x \right\|_r \right), t - pT - \tau \right), \quad \forall t \geq pT + \tau \quad (3.13)$$

Define $\tilde{\sigma}(s, t) := \sigma\left( G^{(p)}(H(s)), t - (p+1)T - r \right) + g\left( \sigma\left( G^{(p)}(H(s)), t - (p+1)T - r - \tau \right) \right)$ for all $t \geq (p+1)T + r + \tau$ and $\tilde{\sigma}(s, t) := G^{(p)}(H(s)) + g\left( G^{(p)}(H(s)) \right)$ for all $t \in [0, (p+1)T + r + \tau)$. Using (3.1), (3.8), (3.10), (3.11) and (3.13) we can conclude that (3.5) holds. Notice that if $\sigma(s, t) = 0$ for all $s \geq 0$ and $t \geq jT$ (i.e., the dead-beat property of order $jT$) then (3.13) implies that $x(t) = 0$ for all $t \geq (j + p + l)T$ (i.e., the dead-beat property of order $(j + p + l)T$). The proof is complete. ◁

Theorem 3.2 can be applied to all forward complete systems satisfying Remark 2.2(d) and Remark 3.1. Here we focus on two special cases: the case of stabilizable Linear Time Invariant (LTI) systems and the case of systems which are Diffeomorphically Equivalent to a Chain of Integrators (DECI). The latter case includes the linearizable strict feedforward systems (see [18]).

**3.2 Design for Stabilizable LTI Systems**

For the Linear Time Invariant (LTI) case, there are matrices $A \in \Re^{n \times n}, B \in \Re^{n \times m}$ such that $f(x, u) = Ax + Bu$. In this case the predictor mapping is given by the explicit expression $\Phi(x, u) = \exp(A(r + \tau))x + \int_{-r-\tau}^{0} \exp(-As)Bu(s)dw$. The linear feedback law $k(x) = Kx$, where $K \in \Re^{m \times n}$, satisfies hypothesis (H3) if and only if there exists $T > 0$ such that all the eigenvalues of the matrix $\exp(AT)\left(I + \int_0^T \exp(-Aw)dw\, BK\right)$ are strictly inside the unit circle on the complex plane.

The use of the sampled-data feedback law $u_i = Ky_i$, where $u_i = u(iT)$ and $y_i = y(iT)$, does not in general guarantee global asymptotic stability for the delayed case, where $r + \tau > 0$. However, Theorem 3.2 can be used for the design of a delay compensator, which guarantees global asymptotic stability for the corresponding closed-loop system. The following corollary is a direct consequence of Theorem 3.2 and its proof is omitted.

**Corollary 3.4 (Stabilization of Linear Networked Control Systems with Delays):** *Let* $T > 0$, $r, \tau \geq 0$ *with* $r + \tau > 0$ *and suppose that there exists* $l \in Z^+$ *such that* $\tau = lT$. *Define* $q := \left[\dfrac{r}{T}\right]$ *and* $\tilde{r} = r - qT$. *Moreover, suppose that all the eigenvalues of the matrix* $\exp(AT)\left(I + \int_0^T \exp(-Aw)dw\, BK\right)$ *are strictly inside the unit circle on the complex plane. Then the closed-loop LTI system*

$$\begin{aligned}\dot{x}(t) &= Ax(t) + Bu(t - \tau) \\ u(t) &= u_i,\ t \in [iT, (i+1)T),\ i \in Z^+\end{aligned} \quad (3.14)$$

*with input applied with zero order hold given by*

$$u_i = K \exp(A(r + \tau))y_i + K \sum_{p=1}^{l+q+1} Q_p B u_{i-p} \quad (3.15)$$

*where* $y_i = y(iT) = x(iT - r)$ *and*

$$\begin{aligned}Q_p &= \exp(ApT)\int_0^T \exp(-As)ds,\ p = 1,\ldots,l+q \\ Q_{l+q+1} &= \exp(A(l+q)T)\left(\int_0^{\tilde{r}} \exp(As)ds\right)\end{aligned} \quad (3.16)$$

*is Globally Exponentially Stable, in the sense that there exist constants* $M, \sigma > 0$ *such that for every* $(x_0, u_0) \in C^0([-r, 0]; \Re^n) \times L^\infty([-r - \tau, 0); \Re^m)$, *the solution* $(x(t), u(t)) \in \Re^n \times \Re^m$ *of system (3.14), (3.15) with initial condition* $\breve{T}_{r+\tau}(0)u = u_0 \in L^\infty([-r - \tau, 0); \Re^m)$, $T_r(0)x = x_0 \in C^0([-r, 0]; \Re^n)$ *satisfies the following inequality for all* $t \geq 0$:

$$\|T_r(t)x\|_r + \|\breve{T}_{r+\tau}(t)u\|_{r+\tau} \leq M\left(\|x_0\|_r + \|u_0\|_{r+\tau}\right)\exp(-\sigma t) \quad (3.17)$$

**Remark 3.5:** If the pair of matrices $(A,B)$ is stabilizable then there exists $K \in \Re^{m \times n}$ such that the matrix $(A+BK)$ is Hurwitz, a symmetric positive definite matrix $P \in \Re^{n \times n}$ and a constant $\mu > 0$ such that $P(A+BK)+(A+BK)'P+\mu I < 0$. Using Corollary 4.3 in [12] with $V(x) = x'Px$, $a(s) := \lambda s$, arbitrary $\lambda \in (0,1)$, and the fact that $A(T,x) \subseteq \left\{ x_0 \in \Re^n : |x-x_0| \leq T|BK||x_0| + T|A|\sqrt{\frac{k_2}{\lambda k_1}}|x| \right\}$, where $k_2 := \max_{|x|=1} x'Px$, $k_1 := \min_{|x|=1} x'Px$, one can show that all the eigenvalues of the matrix

$$\exp(AT)\left(I + \int_0^T \exp(-Aw)dw\, BK\right)$$

are strictly inside the unit circle on the complex plane for all $T > 0$ satisfying

$$T|BK| < 1 \text{ and } 2\frac{T\left(|A|\sqrt{\frac{k_2}{k_1}} + |BK|\right)}{1-T|BK|}|PBK| < \mu$$

Of course, the estimate for the maximum allowable sampling period provided by the above inequalities is conservative in most cases. Other estimates for the maximum allowable sampling period can be found in [31].

**Example 3.6:** We consider the scalar control system

$$\dot{x}(t) = x(t) + u(t) \tag{3.18}$$

where $x(t) \in \Re, u(t) \in \Re$. The system can be exponentially stabilized by the linear feedback $u = -kx$ with $k > 1$ applied with zero order hold, i.e.,

$$\begin{aligned} u(t) &= u_i, \ t \in [iT,(i+1)T) \\ u_i &= -kx(iT), \ i \in Z^+ \end{aligned} \tag{3.19}$$

where the sampling period $T > 0$ must satisfy

$$T < \ln\left(1 + \frac{2}{k-1}\right) \tag{3.20}$$

The use of the same feedback law for the case where measurement delays are present is described by the equations:

$$\begin{aligned} \dot{x}(t) &= x(t) + u(t) \\ u(t) &= u_i, \ t \in [iT,(i+1)T) \\ u_i &= -kx(iT-r), \ i \in Z^+ \end{aligned} \tag{3.21}$$

where $r \geq 0$ is the measurement delay. Numerical experiments for the closed-loop system (3.21) show that for each pair of $k > 1$ and $T > 0$ satisfying (3.20), there exists $r_c > 0$ such that

- if $r < r_c$ then system (3.21) is globally exponentially stable,

- if $r > r_c$ then system (3.21) admits exponentially growing solutions.

For the case $k = 2$, $T = 1$ the value of the critical measurement delay satisfies $r_c \in (0.20, 0.21)$. Figure 1 shows the evolution of the state for system (3.21) with $k = 2$, $T = 1$, $r = 0.1$ and initial condition $x(\theta) = 1$ for $\theta \in [-1, 0]$. The state converges exponentially to zero. Figure 2 shows the evolution of the state for system (3.21) with $k = 2$, $T = 1$, $r = 0.3$ and same initial condition $x(\theta) = 1$ for $\theta \in [-1, 0]$. In this case, the state grows exponentially, indicating instability.

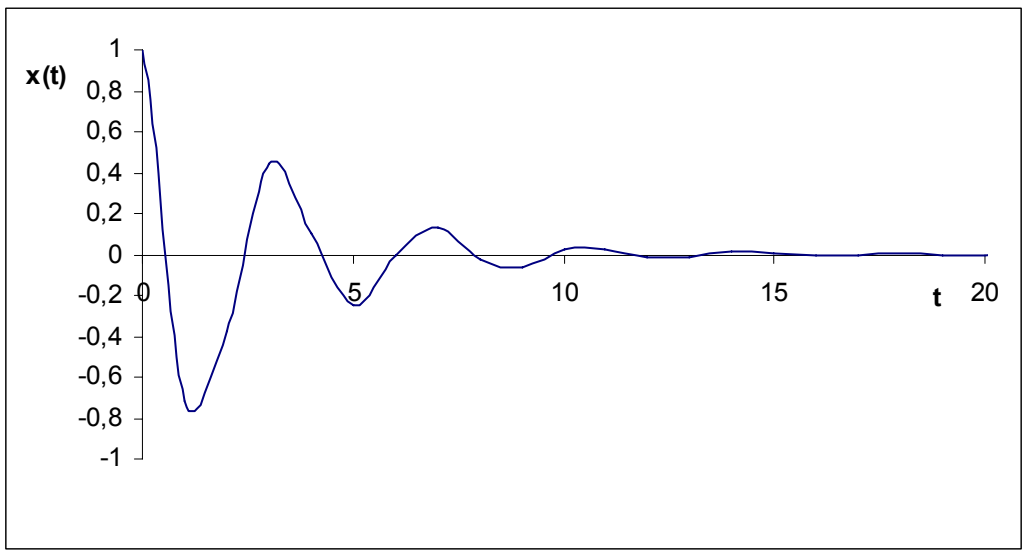

**Figure 1:** The evolution of the state for system (3.21) with $k = 2$, $T = 1$, $r = 0.1$ and initial condition $x(\theta) = 1$ for $\theta \in [-1, 0]$

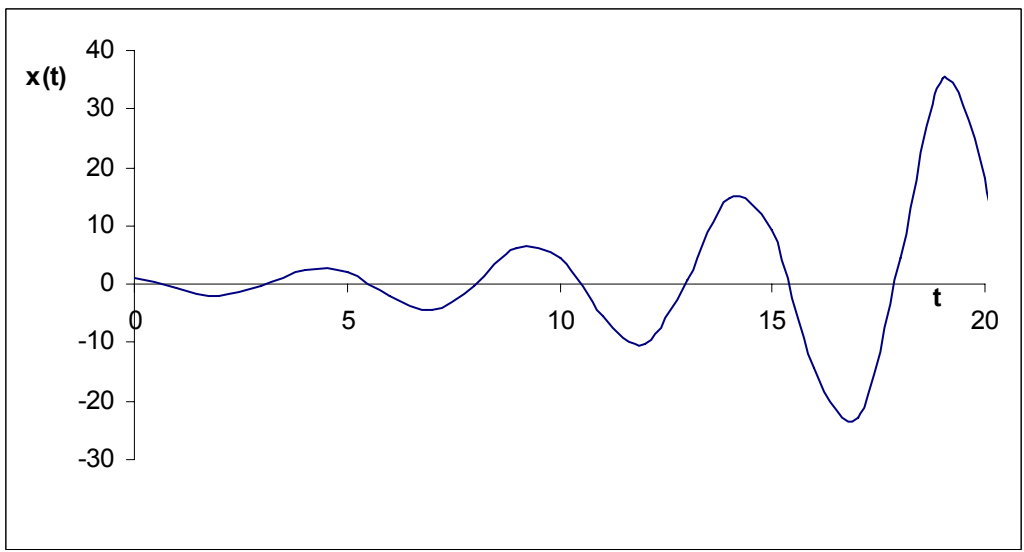

**Figure 2:** The evolution of the state for system (3.21) with $k = 2$, $T = 1$, $r = 0.3$ and initial condition $x(\theta) = 1$ for $\theta \in [-1, 0]$

It is clear that for the case $r > r_c$ one needs a delay compensator. Notice that for the case $k = 2$, $T = 1$ the critical measurement delay $r_c \in (0.20, 0.21)$ is only a small fraction of the sampling period. The usual practice would be to ignore the delay and this would give rise to completely unacceptable results. Corollary 3.4 shows that the feedback law:

$$u(t) = u_i, \; t \in [iT, (i+1)T)$$
$$u_i = -k \exp(r) x(iT - r) - k(\exp(r) - 1) u_{i-1} \qquad (3.22)$$

will guarantee global exponential stability for the closed-loop system (3.18) with (3.22) when $r < T$. Indeed, Figure 3 shows the evolution of the state for the closed-loop system (3.18) with (3.22), $k = 2$, $T = 1$, $r = 0.3$ and initial condition $x(\theta) = 1$ for $\theta \in [-1, 0]$, $u(\theta) = 4$ for $\theta \in [-1, 0)$. The state converges exponentially to zero.

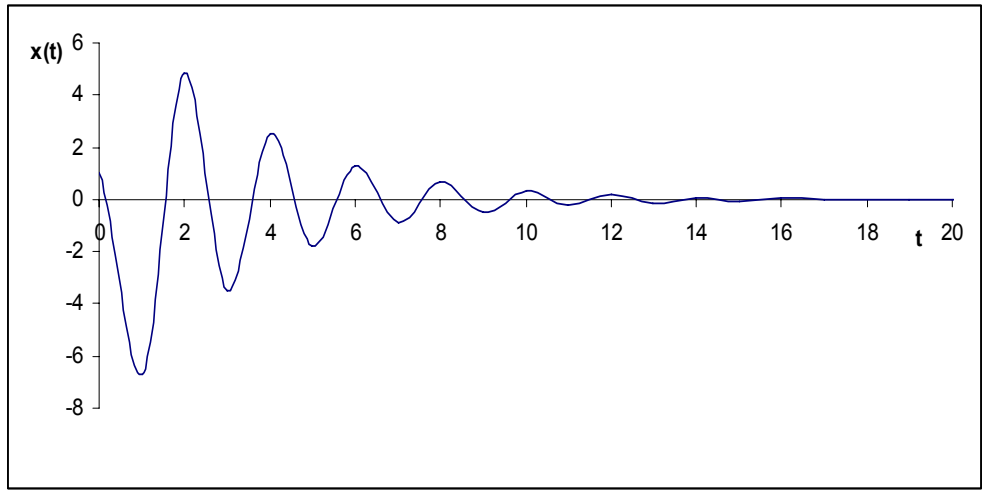

**Figure 3:** The evolution of the state for the closed-loop system (3.18) with (3.22), $k = 2$, $T = 1$, $r = 0.3$ and initial condition $x(\theta) = 1$ for $\theta \in [-1, 0]$, $u(\theta) = 4$ for $\theta \in [-1, 0)$

This example demonstrates that the delay-compensating predictor-based feedback (3.22) extends the range of measurement delays for which stabilization is achieved for given $k > 1$ and $T > 0$ satisfying (3.20).  ◁

## 3.3 Design for Controllable Systems Linearizable by Coordinate Change

The class of systems that are Diffeomorphically Equivalent to a Chain of Integrators (DECI) is the class of all nonlinear systems $\dot{x} = f(x, u)$, $x \in \Re^n, u \in \Re$, where $f : \Re^n \times \Re \to \Re^n$ is locally Lipschitz with $f(0,0) = 0$, for which there exists a global diffeomorphism $\Theta : \Re^n \to \Re^n$ such that the change of coordinates $z = \Theta(x)$ transforms the system to the linear system $\dot{z} = A_0 z + B_0 u$, where $B_0' = (0,...,0,1)$, $A_0 = \{a_{i,j}, i,j = 1,...,n\}$ with $a_{i,i+1} = 1$ for $i = 1,...,n-1$ and $a_{i,j} = 0$ if $j \neq i+1$.

In this case, for every $T > 0$ there exists $K \in \Re^n$, such that all the eigenvalues of the matrix $\exp(A_0 T)\left(I + \int_0^T \exp(-A_0 w) dw\, B_0 K'\right)$ are zero. For example, for $n = 2$ the vector $K \in \Re^2$ is defined by $K' = -\left(\frac{1}{T^2}, \frac{3}{2T}\right)$. If all eigenvalues of the matrix $\exp(A_0 T)\left(I + \int_0^T \exp(-A_0 w) dw\, B_0 K'\right)$ are zero then the sampled-data controller with zero order hold

$$u(t) = K' z(iT),\ t \in [iT, (i+1)T), i \in Z^+$$

applied to the linear system $\dot{z} = A_0 z + B_0 u$ will guarantee the dead-beat property of order $nT$ for the resulting closed-loop system, i.e.,

$$z(t) = 0,\ \text{for all } t \geq nT \text{ and for all initial conditions } z(0) \in \Re^n$$

Thus, we can conclude that the sampled-data controller with zero order hold

$$u(t) = K' \Theta(x(iT)),\ t \in [iT, (i+1)T), i \in Z^+$$

applied to the nonlinear system $\dot{x} = f(x,u)$ will guarantee the dead-beat property of order $nT$ for the resulting closed-loop system.

Therefore, Theorem 3.2 and Corollary 3.4 lead us to the following corollary.

**Corollary 3.7 (Predictor for Linearizable Controllable Systems):** *Let $T > 0$, $r, \tau \geq 0$ with $r + \tau > 0$ and suppose that there exists $l \in Z^+$ such that $\tau = lT$. Define $q := \left[\frac{r}{T}\right]$ and $\tilde{r} = r - qT$. Consider system (1.1) with $m = 1$ and suppose that there exists a global diffeomorphism $\Theta : \Re^n \to \Re^n$ such that*

$$D\Theta(x) f(x,u) = A_0 \Theta(x) + B_0 u,\ \text{for all } x \in \Re^n, u \in \Re \quad (3.23)$$

*where $D\Theta(x)$ is the Jacobian of $\Theta$, $B_0' = (0,...,0,1)$, $A_0 = \{a_{i,j}, i,j = 1,...,n\}$ with $a_{i,i+1} = 1$ for $i = 1,...,n-1$ and $a_{i,j} = 0$ if $j \neq i+1$. Let $K \in \Re^n$ be such that all eigenvalues of the matrix $\exp(A_0 T)\left(I + \int_0^T \exp(-A_0 w) dw\, B_0 K'\right)$ are strictly inside the unit circle on the complex plane. Then the closed-loop system (1.2) with input applied with zero order hold given by*

$$u(t) = u_i,\ t \in [iT, (i+1)T),\ i \in Z^+ \quad (3.24)$$

$$u_i = K' \exp(A_0(r+\tau))\Theta(y_i) + K' \sum_{p=1}^{l+q+1} Q_p B_0 u_{i-p} \tag{3.25}$$

where $y_i = y(iT) = x(iT - r)$ and the matrices $Q_p$ ($p = 1,...,l+q+1$) are defined by (3.16) with $A_0$ in place of $A$, is Globally Asymptotically Stable. Moreover, if all eigenvalues of the matrix $\exp(A_0 T)\left(I + \int_0^T \exp(-A_0 w) dw \, B_0 K'\right)$ are zero then for every $(x_0, u_0) \in C^0([-r,0]; \Re^n) \times L^\infty([-r-\tau,0); \Re^m)$, the solution $(x(t), u(t)) \in \Re^n \times \Re^m$ of system (1.1), (3.24), (3.25) with initial condition $\breve{T}_{r+\tau}(0)u = u_0 \in L^\infty([-r-\tau,0); \Re^m)$, $T_r(0)x = x_0 \in C^0([-r,0]; \Re^n)$ satisfies:

$$x(t) = 0, \text{ for all } t \geq (l+q+1+n)T \tag{3.26}$$

**Example 3.8:** Dead-beat control with a predictor can be applied to any delayed 2-dimensional strict feedforward system, i.e., any system of the form:

$$\dot{x}_1(t) = x_2(t) + p(x_2(t))u(t-\tau), \quad \dot{x}_2(t) = u(t-\tau) \tag{3.27}$$

where $p: \Re \to \Re$ is a smooth function and the measurements are sampled and given by (1.3). The diffeomorphism given by (see [18])

$$\Theta(x) = \left[x_1 - \int_0^{x_2} p(w)dw, \; x_2\right]' \tag{3.28}$$

transforms system (3.27) with $\tau = 0$ to a chain of two integrators. Therefore, the feedback law

$$u = -\frac{1}{T^2}x_1 + \frac{1}{T^2}\int_0^{x_2} p(w)dw - \frac{3}{2T}x_2 \tag{3.29}$$

applied with zero order hold and sampling period $T > 0$ achieves global stabilization of system (3.27) with $\tau = 0$ when no measurement delays are present. Moreover, the dead-beat property of order $2T$ is guaranteed for the corresponding closed-loop system. Interestingly, the feedback law (3.29) is also globally asymptotically stabilizing as a continuous-time controller, placing the closed-loop poles in the $\Theta$-coordinates at $\frac{-3 \pm i\sqrt{7}}{4T}$ for any $T > 0$.

We next consider the case where we have measurement delay $r > 0$ satisfying $r < T$. In this case ($l = q = 0$, $\tilde{r} = r$) we apply Corollary 3.7 and we can conclude that the feedback law (3.24) with

$$u_i = -\frac{1}{T^2}x_1(\tau_i - r) + \frac{1}{T^2}\int_0^{x_2(\tau_i - r)} p(w)dw - \frac{3T+2r}{2T^2}x_2(\tau_i - r) - \frac{r(r+3T)}{2T^2}u_{i-1} \tag{3.30}$$

guarantees the dead-beat property of order $3T$ for the corresponding closed-loop system. Similar formulas to (3.30) are obtained for other cases, where $\tau > 0$ or $r \geq T$. ◁

# 4. Stabilization of a Nonholonomic Mobile Robot Over a Long-Distance Communication Network with Arbitrarily Sparse Sampling

The reduced-order model of a three-wheeled vehicle with two independent rear motorized wheels can be described by the following system of differential equations:

$$\dot{x}(t) = v(t)\cos(\theta(t)), \quad \dot{y}(t) = v(t)\sin(\theta(t)), \quad \dot{\theta}(t) = \omega(t) \tag{4.1}$$

where $(x, y)$ are the coordinates of the center of mass of the vehicle and $\theta$ is the angle between the axis of the vehicle and the horizontal axis. The inputs $v$ and $\omega$ are linear combinations of the angular velocities of the two rear wheels.

The coordinate transformation:

$$x_1 = x\cos(\theta) + y\sin(\theta) \tag{4.2a}$$

$$x_2 = x\sin(\theta) - y\cos(\theta) \tag{4.2b}$$

$$x_3 = \theta \tag{4.2c}$$

and the input transformation:

$$u_1 = v - x_2\omega, \quad u_2 = \omega \tag{4.3}$$

brings system (4.1) to the form (1.4). Many researchers have obtained results for the stabilization of the equilibrium point $0 \in \Re^3$ of system (1.4). Here, we assume that the measurements $x(t-r)$, $y(t-r)$, $\theta(t-r)$, where $r > 0$ ($i = 1,2,3$), are available at discrete time instants which differ by a constant $T > 0$. Moreover, we assume that there is a time delay $\tau \geq 0$ between the computed control action and the applied input (communication delay). In this case the equations of the vehicle are:

$$\dot{x}(t) = v(t-\tau)\cos(\theta(t)), \quad \dot{y}(t) = v(t-\tau)\sin(\theta(t)), \quad \dot{\theta}(t) = \omega(t-\tau) \tag{4.4}$$

with measurements $x(t-r)$, $y(t-r)$ and $\theta(t-r)$.

The reader should notice that hypotheses (H1), (H2) hold for system (4.1). Particularly, there exist smooth time-periodic feedback stabilizers for system (1.4) (see [25,32]) and consequently we can guarantee that hypothesis (H2) holds. The predictor mapping for system (4.4) is given by the following equation:

$$\Phi(x, y, \theta, v, \omega) := \left[ x + \int_{-r-\tau}^{0} v(s)\cos\left(\theta + \int_{-r-\tau}^{s} \omega(p)dp\right)ds \;,\; y + \int_{-r-\tau}^{0} v(s)\sin\left(\theta + \int_{-r-\tau}^{s} \omega(p)dp\right)ds \;,\; \theta + \int_{-r-\tau}^{0} \omega(s)ds \right]'$$

$$\tag{4.5}$$

Using any stabilizing feedback from [25,32] for the nonholonomic integrator (1.4), the coordinate transformation (4.2) and the input transformation (4.3) and Theorem 2.1, we arrive at the following corollary.

**Corollary 4.1:** *Assume that* $k \in C^1(\Re^+ \times \Re^3; \Re^2)$ *with* $k(t,0) = 0$ *for all* $t \geq 0$, *is a time periodic uniform stabilizer for (1.4), i.e., the feedback law* $u_1(t) = k_1(t, x_1(t), x_2(t), x_3(t))$, $u_2(t) = k_2(t, x_1(t), x_2(t), x_3(t))$ *uniformly, globally stabilizes* $0 \in \Re^3$ *for system (1.4). Then for every* $r, \tau \geq 0$, $T > 0$ *the sampled-data dynamic feedback*

$$\dot{z}_x(t) = v(t)\cos(z_\theta(t))$$
$$\dot{z}_y(t) = v(t)\sin(z_\theta(t)), \quad t \in [\tau_i, \tau_{i+1}) \qquad (4.6)$$
$$\dot{z}_\theta(t) = \omega(t)$$

$$z_x(\tau_{i+1}) = x(\tau_{i+1} - r) + \int_{\tau_{i+1}-r-\tau}^{\tau_{i+1}} v(s)\cos\left(\theta(\tau_{i+1} - r) + \int_{\tau_{i+1}-r-\tau}^{s} \omega(p)dp\right) ds$$

$$z_y(\tau_{i+1}) = y(\tau_{i+1} - r) + \int_{\tau_{i+1}-r-\tau}^{\tau_{i+1}} v(s)\sin\left(\theta(\tau_{i+1} - r) + \int_{\tau_{i+1}-r-\tau}^{s} \omega(p)dp\right) ds \qquad (4.7)$$

$$z_\theta(\tau_{i+1}) = \theta(\tau_{i+1} - r) + \int_{\tau_{i+1}-r-\tau}^{\tau_{i+1}} \omega(s) ds$$

*where* $\tau_i = t_0 + iT$, $i \in Z^+$ *and*

$$v(t) = k_1(t + \tau, z_x(t)\cos(z_\theta(t)) + z_y(t)\sin(z_\theta(t)), z_x(t)\sin(z_\theta(t)) - z_y(t)\cos(z_\theta(t)), z_\theta(t)) +$$
$$+ \left(z_x(t)\sin(z_\theta(t)) - z_y(t)\cos(z_\theta(t))\right)\omega(t) \qquad (4.8)$$

$$\omega(t) = k_2(t + \tau, z_x(t)\cos(z_\theta(t)) + z_y(t)\sin(z_\theta(t)), z_x(t)\sin(z_\theta(t)) - z_y(t)\cos(z_\theta(t)), z_\theta(t))$$

*achieves uniform global stabilization of* $0 \in \Re^3$ *for system (4.4).*

At this point it should be emphasized that if the smooth, time-varying feedback proposed in [9] for the stabilization of the nonholonomic integrator were used in (4.8) then the closed-loop system (4.4) with (4.6), (4.7), (4.8) would be non-uniformly in time globally asymptotically stable (see Remark 2.3 above). Moreover, the reader should compare the result of Corollary 4.1 with the results in [17]: no restrictions for the magnitudes of the delays are imposed in the present work.

If the control action is implemented with zero order hold, then a different procedure has to be applied. The reader should notice that for every sampling period $T > 0$ the discontinuous feedback stabilizer

$$k_1(x) = -\frac{2}{T}\text{sgn}(x_2)|x_2|^{1/2} - \frac{2}{T}x_1$$
$$k_2(x) = \frac{1}{T}|x_2|^{1/2}$$
, for $x \in C_1 := \{\xi \in \Re^3 : \xi_2(2\xi_2 - \xi_1\xi_3) \neq 0\}$ (4.9a)

$$k_1(x) = -\frac{x_1 x_3^2}{T(x_1^2 + x_3^2)}$$
$$k_2(x) = \frac{x_3 x_1^2}{T(x_1^2 + x_3^2)}$$
, for $x \in C_2 := \{\xi \in \Re^3 : \xi_2 = 0, \xi_1\xi_3 \neq 0\}$ (4.9b)

$$k_1(x) = -\frac{x_1}{T}$$
$$k_2(x) = -\frac{x_3}{T}$$
, for $x \in C_3 := \{\xi \in \Re^3 : 2\xi_2 = \xi_1\xi_3\}$ (4.9c)

satisfies hypothesis (H3) for system (1.4). To see this notice that inequality (3.1) holds with $g(s) := \frac{2}{T}s + \frac{3}{T}\sqrt{s}$. Furthermore, by explicit computation of the solution one can show that:

- if $x_0 \in C_1$ then the solution $x(t)$ of (1.4) with $u(t) = k(x_0)$ satisfies $x_2(T) = 0$, i.e., $x(T) \in C_2 \cup C_3$,

- if $x_0 \in C_2$ then the solution $x(t)$ of (1.4) with $u(t) = k(x_0)$ satisfies $2x_2(T) = x_1(T)x_3(T)$, i.e., $x(T) \in C_3$,

- if $x_0 \in C_3$ then the solution $x(t)$ of (1.4) with $u(t) = k(x_0)$ satisfies $x(T) = 0$.

It follows that the sampled-data implementation of the feedback (4.9) with sampling period $T > 0$ guarantees the dead-beat property of order $3T$ for the corresponding closed-loop system. The inequality

$$\sup_{0 \leq t < T} |x(t)| \leq 3|x(0)| + T(2 + |x(0)|)\sup_{0 \leq t < T} |u(t)| + \frac{T^2}{2}\sup_{0 \leq t < T} |u(t)|^2$$

which holds for the solution of (1.4) for every initial condition and for every applied input $u \in L^\infty([0,T]; \Re^2)$, in conjunction with inequality (3.1) guarantees for every $T > 0$ the existence of a function $a_T \in K_\infty$ such that for every initial condition the solution of (1.4) with the sampled-data implementation of the feedback (4.9) with sampling period $T > 0$ satisfies $\sup_{0 \leq t < 3T} |x(t)| \leq a_T(|x(0)|)$. The dead-beat property of order $3T$ in conjunction with the previous estimate guarantees that the discontinuous feedback stabilizer defined by (4.9) satisfies hypothesis (H3) for system (1.4). The reader should notice that the feedback design procedure described in [4] can be applied as well for the nonholonomic integrator (1.4) (since (1.4) is an asymptotically controllable homogeneous system with positive minimal power and zero degree).

Therefore, we are in a position to solve the stabilization problem for (4.4) for the case of inputs applied with zero order hold, with measurements $x(t-r)$, $y(t-r)$ and $\theta(t-r)$ available at discrete time instants. Indeed, first we select a sampling period $T>0$ which satisfies $\tau = lT$ for some $l \in Z^+$. Then we can apply Theorem 3.2 for the discontinuous feedback defined by (4.9) and obtain the following result.

**Proposition 4.2:** *Let $r, \tau \geq 0$, $T > 0$ and assume that $\tau = lT$ for some $l \in Z^+$. Then $0 \in \Re^3$ is uniformly globally asymptotically stable for the closed-loop system (4.4) with*

$$\omega(t) = k_2(X\cos(\Theta) + Y\sin(\Theta), X\sin(\Theta) - Y\cos(\Theta), \Theta), \text{ for } t \in [\tau_i, \tau_{i+1}) \quad (4.10)$$

$$v(t) = k_1(X\cos(\Theta) + Y\sin(\Theta), X\sin(\Theta) - Y\cos(\Theta), \Theta) + (X\sin(\Theta) - Y\cos(\Theta))\omega(t), \text{ for } t \in [\tau_i, \tau_{i+1}) \quad (4.11)$$

*where $\tau_i = iT$, $i \in Z^+$, $k: \Re^3 \to \Re^2$ is defined by (4.9) and*

$$\begin{aligned}
X &= x(\tau_i - r) + \int_{\tau_i - r - \tau}^{\tau_i} v(s)\cos\left(\theta(\tau_i - r) + \int_{\tau_i - r - \tau}^{s} \omega(p)dp\right)ds \\
Y &= y(\tau_i - r) + \int_{\tau_i - r - \tau}^{\tau_i} v(s)\sin\left(\theta(\tau_i - r) + \int_{\tau_i - r - \tau}^{s} \omega(p)dp\right)ds \\
\Theta &= \theta(\tau_i - r) + \int_{\tau_i - r - \tau}^{\tau_i} \omega(s)ds
\end{aligned} \quad (4.12)$$

*Moreover, for every initial condition the solution of the closed-loop system (4.4) with (4.10), (4.11), (4.12), where $\tau_i = iT$, $i \in Z^+$, $k: \Re^3 \to \Re^2$ is defined by (4.9) satisfies $x(t) = y(t) = \theta(t) = 0$ for all $t \geq \left(4 + l + \left[\frac{r}{T}\right]\right)T$, where $\left[\frac{r}{T}\right]$ is the integer part of $\frac{r}{T}$.*

It should be noticed that the control action computed by (4.10), (4.11), (4.12) is applied with zero order hold, i.e., $v(t)$ and $\omega(t)$ are constant on each interval $[\tau_i, \tau_{i+1})$, $i \in Z^+$, and hence they are piecewise constant over the integration intervals $[\tau_i - r - \tau, \tau_i]$ in (4.12).

## 5. Concluding Remarks

Stabilization is studied for nonlinear systems with input and measurement delays, and with measurements available only at discrete time instants (sampling times). Two different cases are considered: the case where the input can be continuously adjusted and the case where the input is applied with zero order hold. Under the assumption of forward completeness and certain additional stabilizability assumptions, it is shown that sampled-data feedback laws with a predictor-based delay compensation can guarantee global asymptotic stability for the closed-loop system with no restrictions for the magnitude of the delays. Additionally, when the control is applied continuously and only the measurements are sampled, the sampling time can be arbitrarily long. Applications to the stabilization of linear networked control systems, strict feedforward systems and a nonholonomic mobile robot over a long-distance communication network are presented.

Future work will address the issue of robustness of the proposed feedback laws with respect to actuator and measurement errors, as well as the extension of the obtained results to the case where the delayed and sampled measured output does not necessarily coincide with the state vector.

**Acknowledgments:** The authors would like to thank Professor Costas Kravaris for bringing to their attention the fact that the predictor mapping can be explicitly computed for bilinear systems $\dot{x} = Ax + Bu + uCx$, with $x \in \Re^n, u \in \Re$ and $AC = CA$ (see Remark 2.2(d)).

# Appendix

**Proof of Claim 1 in the proof of Theorem 2.1:** Let $(x_0, z_0, u_0) \in C^0([-r,0]; \Re^n) \times \Re^n \times L^\infty([-r-\tau,0]; \Re^m)$ be arbitrary and consider the solution $(x(t), z(t), u(t)) \in \Re^n \times \Re^n \times \Re^m$ of the closed-loop system (2.5), (2.6), (1.3), (1.1) with initial condition $z(t_0) = z_0 \in \Re^n$, $T_{r+\tau}(t_0)u = u_0 \in L^\infty([-r-\tau,0]; \Re^m)$, $T_r(t_0)x = x_0 \in C^0([-r,0]; \Re^n)$. It is crucial to notice that the solution of $\dot{z}(t) = f(z(t), k(t+\tau, z(t)))$ with $z(t_0) = z_0 \in \Re^n$ satisfies $z(t) = \xi(t+\tau)$, where $\xi(s)$ is the solution of $\dot{\xi}(s) = f(\xi(s), k(s, \xi(s)))$ with $\xi(t_0 + \tau) = z_0 \in \Re^n$. Inequality (2.4) implies that the solution $z(t) \in \Re^n$ exists for all $t \in [t_0, t_0 + T)$ and that the following inequality holds

$$|z(t)| \le \sigma(|z_0|, t - t_0), \quad \forall t \in [t_0, t_0 + T) \tag{A1}$$

It follows that the solution of (2.6) exists for all $t \in [t_0, t_0 + T)$. Continuity of $k \in C^1(\Re^+ \times \Re^n; \Re^m)$ and inequalities (2.3), (A1) imply that the mapping $t \to u(t)$ is continuous on $(t_0, t_0 + T)$ and bounded with $\lim_{t \to t_0^+} u(t) = k(t_0 + \tau, z_0)$ and $\lim_{t \to (t_0+T)^-} u(t) = k(t_0 + T + \tau, z^*)$, where $z^* = \lim_{t \to (t_0+T)^-} z(t)$. Notice that the limit $z^* = \lim_{t \to (t_0+T)^-} z(t)$ exists by virtue of uniform continuity of the mapping $t \to z(t)$ on $[t_0, t_0 + T)$. By virtue of inequalities (2.3) and (A1) we obtain the inequality

$$\sup_{t_0 - \tau - r \le s < t_0 + T} |u(s)| \le \|T_{r+\tau}(t_0)u\|_{r+\tau} + g(\sigma(|z_0|, 0)) \tag{A2}$$

Using hypothesis (H1) we may conclude that the solution $x(t) \in \Re^n$ of (1.1) exists for all $t \in [t_0, t_0 + T]$. Indeed, by virtue of the Fact, we can guarantee the existence of $\zeta \in K_\infty$ such that

$$|x(t)| \leq \zeta\left(|x(t_0)| + \sup_{t_0-\tau \leq s < t-\tau}|u(s)|\right), \text{ for all } t \in (t_0, t_0+T]$$

The above inequality in conjunction with (A2) and the trivial inequality $\|T_r(t)x\|_r \leq \|T_r(t_0)x\|_r + \sup_{t_0 \leq s \leq t}|x(s)|$ gives:

$$\|T_r(t)x\|_r \leq \|x_0\|_r + \zeta\left(\|x_0\|_r + \|T_{r+\tau}(t_0)u\|_{r+\tau} + g(\sigma(|z_0|,0))\right), \forall t \in [t_0, t_0+T] \quad (A3)$$

Finally, we define $u(t_0+T) = k(t_0+T+\tau, z(t_0+T))$, where $z(t_0+T) = \Phi(x(t_0+T-r), \breve{T}_{r+\tau}(t_0+T)u)$. Again, using (2.1), (2.3), (A2) and (A3) we obtain:

$$|z(t_0+T)| \leq a\left(\|T_{r+\tau}(t_0)u\|_{r+\tau} + g(\sigma(|z_0|,0)) + \|x_0\|_r + \zeta(\|x_0\|_r + \|T_{r+\tau}(t_0)u\|_{r+\tau} + g(\sigma(|z_0|,0)))\right)$$

and

$$|u(t_0+T)| \leq g\left(a\left(\|T_{r+\tau}(t_0)u\|_{r+\tau} + g(\sigma(|z_0|,0)) + \|x_0\|_r + \zeta(\|x_0\|_r + \|T_{r+\tau}(t_0)u\|_{r+\tau} + g(\sigma(|z_0|,0)))\right)\right)$$

Using (A1), (A2), (A3) and the above inequalities we are in a position to construct a function $G \in K_\infty$ such that inequality (2.8) holds. The proof is complete. ◁

**Proof of Claim 3 in the proof of Theorem 3.2:** Consider first the case $\tau > 0$. Since $\tau = lT$ for some $l \in Z^+$, it follows that $\tau \geq T$. By virtue of the Fact, we can guarantee the existence of $b \in K_\infty$ such that

$$|x(t)| \leq b\left(|x(\tau)| + \|\breve{T}_{r+\tau}(\tau)u\|_{r+\tau}\right), \text{ for all } t \in [\tau, \tau+T] \quad (A3)$$

Moreover, using the equation $u(t) = k(\Phi(x(\tau-r), \breve{T}_{r+\tau}(\tau)u)), t \in [\tau, \tau+T)$, inequalities (2.1) and (3.1), we obtain

$$|u(t)| \leq g\left(a\left(\|T_r(\tau)x\|_r + \|\breve{T}_{r+\tau}(\tau)u\|_{r+\tau}\right)\right), \text{ for all } t \in [\tau, T+\tau) \quad (A4)$$

Finally, using the trivial inequalities $\|\breve{T}_{r+\tau}(t)u\|_{r+\tau} \leq \|\breve{T}_{r+\tau}(\tau)u\|_{r+\tau} + \sup_{\tau \leq s < t}|u(s)|$ and $\|T_r(t)x\|_r \leq \|T_r(\tau)x\|_r + \sup_{\tau \leq s \leq t}|x(s)|$ in conjunction with (A3), (A4), we can conclude that (3.9) holds with $G(s) := s + b(s) + g(a(s))$.

Next consider the case $\tau = 0$. Using the equation $u(t) = k(\Phi(x(-r), \breve{T}_r(0)u)), t \in [0, T)$, inequalities (2.1) and (3.1), we obtain (A4) with $\tau = 0$. By virtue of the Fact and the assumption $\tau = 0$, there exists $b \in K_\infty$ such that

$$|x(t)| \leq b\left(|x(0)| + \sup_{0 \leq s < t}|u(s)|\right), \text{ for all } t \in (0, T] \quad (A5)$$

Combining (A4) and (A5) we obtain:

$$|x(t)| \leq b\left(|x(0)| + g\left(a\left(\|T_r(0)x\|_r + \|\breve{T}_r(0)u\|_r\right)\right)\right), \text{ for all } t \in [0,T] \quad (A6)$$

Again, using the trivial inequalities $\left\|\tilde{T}_r(t)u\right\|_r \leq \left\|\tilde{T}_r(0)u\right\|_{r+\tau} + \sup_{0\leq s<t}|u(s)|$ and $\left\|T_r(t)x\right\|_r \leq \left\|T_r(0)x\right\|_r + \sup_{0\leq s\leq t}|x(s)|$ in conjunction with (A6) and (A4) with $\tau = 0$, we can conclude that (3.9) holds with $G(s) := s + b(s + g(a(s))) + g(a(s))$. The proof is complete. ◁